\documentclass[preprint,11pt,3p]{elsarticle}

\usepackage{amssymb}
\usepackage{amsthm}
\usepackage{amsmath,  bbm}

\usepackage{comment}

\usepackage{hyperref}

\newcommand{\rd}{{\mathbb R}^d}
\newcommand{\zd}{{\mathbb Z}^d}
\newcommand{\td}{{\mathbb T}^d}
\newcommand{\n}{{\mathbb N}}

\newcommand{\h} {\widehat}
\newcommand{\nul} {{\bf 0}}
\newtheorem{theo}{Theorem}
\newtheorem{lem}[theo]{Lemma}

\newtheorem {defi} [theo] {Definition}
\newtheorem {rem} [theo] {Remark}

\journal{*}

\begin{document}

\begin{frontmatter}

\title{Directional Heisenberg uncertainty product \tnoteref{label0}}
\tnotetext[label0]{The work is supported by Volkswagen Foundation}

\author[label1]{A. Krivoshein\corref{cor1}}
\address[label1]{St. Petersburg State University}
\address[label3]{Donetsk National University}
\address[label4]{University of L\"ubeck}

\cortext[cor1]{Corresponding author}

\ead{krivosheinav@gmail.com}
%
\author[label1]{E. Lebedeva}
\ead{ealebedeva2004@gmail.com}%

\author[label3]{E. Neiman}
\ead{evg\_sqrt@mail.ru}
\author[label4]{J. Prestin}
\ead{prestin@math.uni-luebeck.de}%

\begin{abstract}
A directional time-frequency localization measure for functions defined on the $d$-dimensional Euclidean space is introduced. A connection between this measure and its periodic counterpart is established. For a class of functions, an optimization problem for finding the optimal direction, along which a function is best or worst localized, is solved.  
\end{abstract}


\end{frontmatter}


\section{Introduction}
\label{sec1}

The paper continues the investigation of the properties of the directional uncertainty product, that was recently introduced for the periodic case in~\cite{KLParhiv}. This paper deals with a non-periodic counterpart. In the framework of the standard operator approach (see, e.g., Selig in~\cite{Selig2002UC} or Goh, Micchelli in~\cite{GohMicch2002UC}) we introduce a pair of operators, that are appropriate for measuring a time-frequency localization along directions for functions defined on $\rd$. The corresponding uncertainty principle is valid automatically, the lower bound of the directional uncertainty product is equal to $1/4$ and is attained on the class of functions, that are Gaussian exponentials up to a multiplication on arbitrary smooth functions. 
Our definition, in contrast to definitions given by Goh and Goodman in~\cite{GohGoodman2004UC}, Ozawa and Yuasa in~\cite{OzawaYuasa16UC}, includes the directionality explicitly in a natural way.

We establish a connection between the directional uncertainty products in the periodic and non-periodic case (see Subsection 3.1). Namely, for an appropriate class of functions $f$, the periodic directional uncertainty product of its periodization  tends to the non-periodic directional uncertainty product of $f$ as the period goes to infinity. This connection is also established for the uncertainty product, that was suggested by Goh and Goodman  in~\cite{GohGoodman2004UC}. We also study the dependence on the direction of the directional uncertainty product for a fixed function (see Subsection~3.2).  It is an optimization problem, one needs to find a direction along which the directional uncertainty product has its minimum or maximum. For a class of symmetric functions the optimization problem is solved analytically. Finally, by using the Fourier-Hermite series, we state for a class of symmetric functions that the lower bound of the directional uncertainty product can be improved (see Subsection 3.3). 
The proofs of all statements are given in Section 4. Several examples illustrating the results of Subsection 3.2 are placed in Section 5. 

\section{Basic notations and definitions}
\label{sec2}

We use the standard multi-index notations.
    Let $d\in\n,$ $\rd$ be the $d$-dimensional Euclidean space, $\{e_j, 1\le j\le d\}$ be the standard basis in $\rd$, $\zd$ is the integer lattice  in $\rd$, 
    $\td=\rd\slash\zd$ be the $d$-dimensional torus.
    Let  $x = (x_1,\dots, x_d)^{\mathrm T}$ and 
    $y =(y_1,\dots, y_d)^{\mathrm T}$ be column vectors in $\rd$.
    Then $\langle x, y\rangle:=x_1y_1+\dots+x_dy_d$,
    $\|x\| := \sqrt {\langle x, x\rangle}$.  
    We say that
     $x \geq y,$ if $x_j \geq y_j$ for all $j= 1,\dots, d,$
    and we say that 
 	$x > y,$ if $x \geq y $ and $x\neq y.$
    $\zd_+:=\{\alpha\in\zd:~\alpha\geq~{\bf0}\},$
    where  ${\bf0}=(0,\dots, 0)$ denotes the origin in $\rd$.
    For $\alpha=(\alpha_1,\dots,\alpha_d)^{\mathrm T}\in \zd_+$,
    denote $|\alpha|:=\alpha_1+ \dots +\alpha_d.$ $\mathbbm{1}_{K} (x)$ is the characteristic (indicator) function of a set $K\subset \rd$.
    
	For a smooth enough function $f$ defined on $\rd$ and a multi-index $\alpha\in\zd_+$, 
	$D^{\alpha} f$ denotes the derivative of $f$ of order $\alpha$ and
	$D^{\alpha} f =\frac{\partial^{|\alpha|} f}
    {\partial x^{\alpha}}=\frac{\partial^{|\alpha|} f
    }{\partial^{\alpha_1}x_1\dots
    \partial^{\alpha_d}x_d}.$
	The directional derivative of a smooth enough function $f$ defined on $\rd$ along a vector $L=(L_1,...,L_d)\in\rd \setminus \{\nul\}$ is denoted by
	$\frac{\partial f}{\partial L} =\sum_{j=1}^d L_j \frac{\partial f}{\partial x_j}.
	$ 

For a function $f\in L_2(\td)$ its norm is denoted by 
    $\|f\|^2_{\td} = \int_{\td} |f(x)|^2 \mathrm{d}x$.
    The Fourier coefficients of a function $f\in L_2(\td)$ are given by $c_k = c_k(f) = \int_{\td} f(x) \mathrm{e}^{-2\pi \mathrm{i} \langle k,x\rangle} \mathrm{d}x$, $k\in\zd.$
    For a function $f\in L_2(\rd)$ its norm is denoted by $\|f\|_2^2 = \int_{\rd} |f(x)|^2 \mathrm{d}x.$
    The Fourier transform of a function $f\in L_1(\rd)\bigcap L_2(\rd)$ is given by $\h f(\xi) = \int_{\rd} f(x) \mathrm{e}^{-2\pi \mathrm{i} \langle x,\xi \rangle} \mathrm{d}x$ and can be naturally extended to $L_2(\rd)$.
    The Sobolev space  $H^1(\rd)$ consists of functions in $L_2(\rd)$ such that all its derivatives of the first order are also in $L_2(\rd)$. Analogously we define $H^1(\td)$. Note that
    $$
    H^1(\rd)= \left\{ f\in L_2(\rd): \int_{\rd} \|\xi\|^2 |\h f(\xi)|^2 \mathrm{d}\xi < \infty \right\}.
    $$

Let ${\cal H}$ be a Hilbert space with inner product $\langle \cdot, \cdot \rangle$ and with norm $\|\cdot\| := \langle \cdot, \cdot \rangle^{1/2}$. Let ${\cal A}$, ${\cal B}$ be two linear operators with domains ${\cal D}({\cal A})$, ${\cal D}({\cal B})\subseteq {\cal H}$ and ranges in ${\cal H}$. The variance of non-zero $f\in {\cal D}({\cal A})$ with respect to the operator ${\cal A}$ is defined to be
 $$\Delta_{}({\cal A},f) = \|{\cal A} f\|^2 - \frac{|\langle {\cal A} f, f\rangle|^2}{\|f\|^2}.$$ 
The commutator of ${\cal A}$ and ${\cal B}$ is defined  by 
$
[{\cal A},{\cal B}] := {\cal A}{\cal B}-{\cal B}{\cal A}
$
with domain $ {\cal D}({\cal A}{\cal B}) \bigcap {\cal D}({\cal B}{\cal A})$.

An operator approach for the definition of the uncertainty principle for self-adjoint operators was established by Folland in~\cite{Folland1989Book}. This approach was extended to two normal or symmetric operators by Selig in~\cite{Selig2002UC}  and Goh, Micchelli in~\cite{GohMicch2002UC}. For several operators this approach was generalized by Goh and Goodman in~\cite{GohGoodman2004UC}.
\begin{theo}~\cite[Theorem 4.1]{GohGoodman2004UC}
    Let ${\cal A}_1,\dots{\cal A}_n,$ ${\cal B}_1,\dots{\cal B}_n$ be symmetric or normal operators with domain and range in the same Hilbert space ${\cal H}$. Then for any non-zero $f$ in ${\cal D}({\cal A}_j{\cal B}_j) \bigcap {\cal D}({\cal B}_j{\cal A}_j)$, $j=1,\dots,n$,
    $$
    \frac{1}{4} \left( \sum_{j=1}^n |\langle [{\cal A}_j,{\cal B}_j]  f, f\rangle| \right)^2 \le \left(\sum_{j=1}^n \Delta_{}({\cal A}_j,f)\right)\left(\sum_{j=1}^n \Delta_{}({\cal B}_j,f)\right).
    $$
\end{theo}

If the commutator $\langle [{\cal A}_j,{\cal B}_j]  f, f\rangle$  is non-zero for all $j=1,\dots,n$, then the uncertainty product for $f$ is defined as
$$\mathrm{UP}^{{\cal H}}(f):= \left(\sum_{j=1}^n \Delta_{}({\cal A}_j,f)\right)\left(\sum_{j=1}^n \Delta_{}({\cal B}_j,f)\right) \left( \sum_{j=1}^n |\langle [{\cal A}_j,{\cal B}_j]  f, f\rangle| \right)^{-2} .$$

The well-known Heisenberg uncertainty product for functions in $L_2({\mathbb R})$ fits in this operator approach, if $n=1$ and the two operators are as follows ${\cal A} f(x)= 2 \pi x f(x)$, ${\cal B} f(x)= \frac{\mathrm{i} }{2\pi}\frac{\mathrm{d}f}{\mathrm{d}x}(x)$. Their commutator is $[{\cal A},{\cal B}]=-\mathrm{i}{\cal I}$, where ${\cal I}$ is the identity operator. Both operators are self-adjoint on their domains.
The Heisenberg uncertainty product characterizes the time-frequency localization of a function and the uncertainty principle states that any function cannot have arbitrary good localization in both time and frequency domain. 
It is known that the Heisenberg uncertainty product attains its minimum when $f$ is the Gaussian function.

The Breitenberger uncertainty product is defined for the space of periodic functions $L_2({\mathbb T})$. In this case, ${\cal A}^{{\mathbb T}} f(x)= \mathrm{e}^{2\pi \mathrm{i} x} f(x)$, ${\cal B}^{{\mathbb T}} f(x)= \frac{ \mathrm{i} }{2\pi}\frac{\mathrm{d}f}{\mathrm{d}x}(x)$.

There were several attempts to define the uncertainty product for the multivariate periodic and non-periodic cases. For instance, Goh and Goodman in~\cite{GohGoodman2004UC} suggested to take a collection of operators, where each operator is responsible for one coordinate (or variable). 
For the non-periodic case, these operators are
${\cal A}_j f(x) = 2 \pi x_j f(x)$, ${\cal B}_j f(x)= \frac{ \mathrm{i} }{2\pi} \frac{\partial f}{\partial x_j}(x)$, the commutator is  $[{\cal A}_j,{\cal B}_j]=-\mathrm{i}{\cal I}$ $j=1,\dots,d$, $x\in\rd.$ The corresponding uncertainty product is defined as
\begin{equation}
  \mathrm{UP}_{\textrm{GG}}(f):=\frac{1}{d^2 \|f\|^4_2} \sum_{j=1}^d \Delta_{}({\cal A}_j,f)\sum_{k=1}^d \Delta_{}({\cal B}_k,f)
  \label{fUPGG}
\end{equation}
and it attains its minimum at the multivariate Gaussian function $f(x) = a \mathrm{e}^{-\|b x- c\|^2}$, $a,b \in \rd \setminus\{\nul\},$ $c\in\rd.$
Also, some other approaches were suggested by 
Ozawa and Yuasa in~\cite{OzawaYuasa16UC}.

In fact, the above approaches for the definition of the uncertainty product do not deal with a new phenomenon, that appears in the multidimensional case, namely, the localization of a function along a particular direction. We suggest an approach that allows to include this directionality into the definition.

The directional uncertainty product for $\rd$ along a direction $L\in\rd$ we define using two operators
$$
{\cal A}_L f(x) =   2 \pi \langle L, x\rangle f(x), \quad 
{\cal B}_L f(x) = \frac{ \mathrm{i} }{2\pi} \frac{\partial f}{\partial L} (x).
$$
Note that the domains of these operators are
$$
{\cal D}({\cal A}_L) = \{ f\in L_2(\rd): \int_{\rd} \|x\|^2 |f(x)|^2 \mathrm{d}x < \infty\} = \{ f\in L_2(\rd): \h f \in H^1(\rd)\},$$
and ${\cal D}({\cal B}_L) = H^1(\rd).$
Both operators are self-adjoint.
The commutator is  $[{\cal A}_L,{\cal B}_L]= - \mathrm{i}  \|L\|^2 {\cal I}.$ Hence,  for any non-zero $f\in {\cal D}({\cal A}_L{\cal B}_L) \bigcap {\cal D}({\cal B}_L{\cal A}_L)$
\begin{equation}
    \mathrm{UP}_L (f) := \frac{\Delta_{}({\cal A}_L,f) \Delta_{}({\cal B}_L,f)}{\|L\|^4  \|f\|^4_{2}} \geq \frac 14.
    \label{fUPRD}
\end{equation}
The uncertainty principle is valid automatically, due to the operator approach. Clearly, $\mathrm{UP}_L(f)$ is well-defined for the wider class of functions  $f\in {\cal D}({\cal A}_L) \bigcap {\cal D}({\cal B}_L)$ and by density arguments  also $ \mathrm{UP}_L (f)\geq \frac 14$,  since the variances are continuous functionals on their domains.

The main purpose of this paper is to study the properties of the directional uncertainty product.   


\section{Properties of the directional uncertainty product}

 First of all, we note that modifications of a function like shifts, modulations, scaling and replacing the function by its Fourier transform do not change $\mathrm{UP}_L.$ The directional uncertainty product of a rotated function is equal to the uncertainty product of the initial function along a rotated directional vector.   
\begin{lem}
    Let $f\in {\cal D}({\cal A}_L) \bigcap {\cal D}({\cal B}_L)$.
    Then
    \begin{enumerate}
        \item 
        if $g(x)=a \mathrm{e}^{2\pi \mathrm{i}\langle W, x\rangle} f(b x-x_0),$ where $a,b\in {\mathbb R}$, $x_0, W\in\rd,$ or $g = \h f,$ then 
    $\mathrm{UP}_L(g) = \mathrm{UP}_L(f).$
    
    \item if $U\in{\mathbb R}^{d\times d}$ is a unitary matrix and $g(x) := f(Ux)$ then
    $$
    \Delta({\cal A}_L, g) = \Delta({\cal A}_{UL} f),\ \ \ 
    \Delta({\cal B}_L, g) = \Delta({\cal B}_{UL}, f).
    $$
    \end{enumerate}
    \label{lem:UCRdSTD}
\end{lem}

The proof can be done by straightforward computations.

Next, we establish the set of optimal functions for $\mathrm{UP}_L,$ i.e. $f$ such that $\mathrm{UP}_L(f)=\frac{1}{4}.$


\begin{lem}
Let $L\in\rd,$ $\|L\|=1$, $\mu\in{\mathbb R}\setminus\{0\}.$ For a function $f$ defined by
$$
f(x) = \mathrm{e}^{-\frac{2 \pi^2 }{\mu}\langle L,x\rangle^2 } \Phi (L_2 x_1 - L_1 x_2,L_3 x_1 - L_1 x_3,\dots,L_d x_1 - L_1 x_d),
$$
where $\Phi$ is an arbitrary continuously differentiable function (such that $\mathrm{UP}_L(f)$ makes sense),
it is valid that $\mathrm{UP}_L(f) = \frac{1}{4}.$
\label{lemExactUP}
\end{lem}

\subsection{Connection between periodic and non-periodic case}

In this subsection, we establish a connection between the directional uncertainty products in periodic and non-periodic cases. In the univariate case, this connection between the Heisenberg and Breitenberger uncertainty products was stated  
in~\cite{PrestinQRS2003UC}.

The counterpart of the directional uncertainty product for the periodic case was introduced in~\cite{KLParhiv}. 
It is defined using the operators  
  $${\cal A}_L^{\td} f(x) = \mathrm{e}^{2\pi \mathrm{i} \langle L, x\rangle} f(x), \quad {\cal B}_L^{\td} f(x) = \frac{ \mathrm{i} }{2\pi} \frac{\partial f}{\partial L}(x),$$
 $L\in\zd \setminus\{\nul\}$.
The domains of these operators are
${\cal D}({\cal A}^{\td}_L) = L_2(\td),$ ${\cal D}({\cal B}^{\td}_L) = H^1(\td)$
and  ${\cal A}^{\td}_L$ is normal, ${\cal B}^{\td}_L$ is self-adjoint. 
The commutator for $f\in H^1(\td)$ is $      [{\cal A}^{\td}_L,{\cal B}^{\td}_L]f= \|L\|^2 {\cal A}^{\td}_L f.$
Thus, the directional uncertainty product for a function $f\in  H^1(\td)$ such that $ {\cal A}^{\td}_L f \neq 0$ is defined as 
    {\small  \begin{equation}
  \mathrm{UP}^{\td}_L(f) = \frac{1}{\|L\|^4}\left( \frac{\|f\|^4_{\td}}{|\langle{\cal A}^{\td}_Lf,f\rangle_{\td}|^2} - 1 \right)\left( \frac{\|{\cal B}^{\td}_L f\|^2_{\td}}{\|f\|^2_{\td}} - \frac{|\langle{\cal B}^{\td}_L f, f\rangle_{\td}|^2}{\|f\|^4_{\td}} \right):= \frac{1}{\|L\|^4}\text{var}_L^{\mathrm A}(f) \text{var}_L^{\mathrm F}(f),
  \label{fUPTD}
  \end{equation} }
  where $\text{var}_L^{\mathrm A}(f)$ is the angular directional variance and $\text{var}_L^{\mathrm F}(f)$ is the frequency directional variance.
Also, we introduce the notion of admissible functions, for which the connection will be valid. 

\begin{defi}
A non-zero function $f\in L_2(\rd)$ is called admissible if $f$ is continuously differentiable up to order one, $f\in H^1(\rd)$, $\h f\in H^1(\rd)$ and
\begin{eqnarray*}
|f(x)| &\le & \frac{C_1}{\|x\|^\gamma}, \quad \text{ for all } x\in\rd,
\\
\left| \frac{\partial f}{\partial x_j}(x)\right| &\le & \frac{C_2}{\|x\|^\beta}, \quad \text{ for all } x\in\rd,\ \  j=1,\dots,d,
\end{eqnarray*}
where $C_1>0$ and $C_2>0$ are some constants, $\beta>d,$ $\gamma > \max\{\frac{d}{2}+1,d\}$. 
\end{defi}

For an admissible function $f$ and a parameter $\lambda \in {\mathbb R}$, we denote $f_\lambda(x) : = \sqrt{\lambda^d} f(\lambda x).$ The function $f_\lambda$ is also admissible. 
Consider the periodized version of a scaled admissible function, namely
$$
f_\lambda^{\mathrm{per}} (x) := \sqrt{\lambda^d} \sum_{k\in\zd} f(\lambda(x+k)) =  \sum_{k\in\zd} f_\lambda(x+k). 
$$
In Section 4 it is proved that $f_\lambda^{\mathrm{per}}$ and $\frac{\partial f_\lambda^{\mathrm{per}}}{\partial L}$ are continuous functions in $L_2(\td).$ For admissible functions we can state a connection between $ \mathrm{UP}_L^{\td}$ and $ \mathrm{UP}_L.$

\begin{theo}
    Let $f$ be admissible, $L\in\zd\setminus\{\nul\}$ and $\lambda>0.$ Then
    $$
    \lim_{\lambda\to\infty} \mathrm{UP}_L^{\td} (f_\lambda^{\mathrm{per}}) = \mathrm{UP}_L(f).
    $$
    \label{theoL}
\end{theo}

For the space $L_2(\td)$ of multivariate periodic functions, Goh and Goodman in~\cite{GohGoodman2004UC} suggested to take the operators as follows ${\cal A}^{\td}_j f(x) = \mathrm{e}^{2\pi \mathrm{i} x_j} f(x)$, ${\cal B}^{\td}_j f(x)= \frac{ \mathrm{i} }{2\pi} \frac{\partial f}{\partial x_j}(x)$, $j=1,\dots,d.$ Note that the domains of these operators are $\bigcap_{j=1}^d {\cal D}({\cal A}^{\td}_j) = L_2(\td),$ $\bigcap_{j=1}^d {\cal D}({\cal B}^{\td}_j) = H^1(\td)$. Operators ${\cal A}^{\td}_j$ are normal, ${\cal B}^{\td}_j$ are self-adjoint.
The commutators for $f\in H^1(\td)$ are $[{\cal A}^{\td}_j,{\cal B}^{\td}_j]f=  {\cal A}^{\td}_j f.$
If the commutator $\langle [{\cal A}^{\td}_j,{\cal B}^{\td}_j]  f, f\rangle$  is non-zero for all $j=1,\dots,d$, then the uncertainty product for $f$ is defined as
\begin{multline}
    \mathrm{UP_{GG}^{\td}}(f):= \left(\sum_{j=1}^d \Delta_{}({\cal A}^{\td}_j,f)\right)\left(\sum_{j=1}^d \Delta_{}({\cal B}^{\td}_j,f)\right) \left( \sum_{j=1}^d |\langle [{\cal A}^{\td}_j,{\cal B}^{\td}_j]  f, f\rangle| \right)^{-2} 
\\
 =  {\small  
\frac{ \sum\limits_{j=1}^d \left(\|f\|^4_{\td} - |\langle{\cal A}_j^{\td} f,f\rangle|^2\right)} {\left( \sum\limits_{j=1}^d |\langle{\cal A}_j^{\td} f,f\rangle|\right)^2} \sum\limits_{j=1}^d \left( \frac{\|{\cal B}_j^{\td} f\|^2_{\td}}{\|f\|^2_{\td}} - \frac{|\langle{\cal B}_j^{\td} f, f\rangle|^2}{\|f\|^4_{\td}} \right):= \text{var}_{\mathrm{GG}}^{\mathrm A}(f) \text{var}_{\mathrm{GG}}^{\mathrm F}(f).
  }
\label{fUPGGTD}
\end{multline}

In these terms, the uncertainty principle says that the uncertainty product $\mathrm{UP_{GG}^{\td}}(f)$ cannot be smaller than $\frac{1}{4}$ for any appropriate function $f$.
In this case, the connection between $\mathrm{UP_{GG}^{\td}}$ and $\mathrm{UP_{GG}}$ is also valid.
    
\begin{theo}
    Let $f$ be admissible, $L\in\zd$ and $\lambda>0.$ Then
    $$
    \lim_{\lambda\to\infty} \mathrm{UP_{GG}^{\td}} (f_\lambda^{\mathrm{per}}) = \mathrm{UP_{GG}}(f).
    $$
    \label{theoGG}
\end{theo}

\subsection{Dependence of a localization on the direction for a fixed function}
In this subsection, we fix a function  $f\in {\cal D}({\cal A}_L) \bigcap {\cal D}({\cal B}_L)$ and study how the uncertainty product of this function  depends on  a direction $L\in\rd$. 
Denote
$$
\alpha_{L}(f):=\frac{\langle {\cal A}_L f,f\rangle}{\|f\|_2^2},\ \ \ 
\beta_{L}(f):=\frac{\langle {\cal B}_L f,f\rangle}{\|f\|_2^2},
$$
so time and frequency variances take the form 
$$
\Delta({\cal A}_L, f) = \|{\cal A}_L f\|_2^2 - |\alpha_{L}(f)|^2 \|f\|_2^2,
\quad
\Delta({\cal B}_L, f) = \|{\cal B}_L f\|_2^2 - |\beta_{L}(f)|^2 \|f\|_2^2. 
$$
Without loss of generality we set $\|f\|_2=1$ and $\|L\|=1.$
In the next theorem we give a complete analytic solution for the following extremal problems   
${\rm min}_{\|L\|=1} \mathrm{UP}_L (f)$ and 
${\rm max}_{\|L\|=1} \mathrm{UP}_L (f)$, as the function $f$ satisfies a special type of symmetry relations (see formulas (\ref{symm}) below).

\begin{theo}
\label{1depend_L}
Let $f\in {\cal D}({\cal A}_L) \bigcap {\cal D}({\cal B}_L)$, $\|f\|_2=1$, and 
\begin{gather}
\notag
|f(x_1, \dots, x_k, \dots, x_n)|=|f(x_1, \dots, -x_k, \dots, x_d)|, \\
|\widehat{f}(x_1, \dots, x_k, \dots, x_n)|=|\widehat{f}(x_1, \dots, -x_k, \dots, x_d)|
\label{symm}
\end{gather}
 for all  $k=1,\dots,d.$
Denote 
$$
M_k= (2\pi)^2  \int_{\rd}x_k^2 |f(x)|^2\,\mathrm{d}x  \ \ \mbox{ and } \ \   \widehat{M}_k=\int_{\rd}x_k^2 |\widehat{f}(x)|^2\,\mathrm{d}x. 
$$
Let $A$ be a $d\times d$ matrix whose elements are $(M_k\widehat{M_j}+M_j\widehat{M_k})/2$, $j, k =1,\dots,d$. Let $A_{j_1,\dots,j_q}$ be a submatrix, cut down from $A$ by removing its $j_1$-th, $\dots$, $j_q$-th  row and $j_1$-th, $\dots$, $j_q$-th column, $q=1, \dots, d-1$. Denote ${\bf A}$ the set of all those matrices $A_{j_1,\dots,j_q}$ whose determinant is  not equal to zero, and all the coordinates of the vector $A_{j_1,\dots,j_q}^{-1}E$ are nonnegative, $E=(1,\dots,1)\in{\mathbb R}^{d-q}.$   
 Then 
$$
\min_{\|L\|=1} \mathrm{UP}_L (f) = \min_{L\in {\cal L}} \mathrm{UP}_L (f)    \ \ \ 
\mbox{ and }
\ \ \ 
\max_{\|L\|=1} \mathrm{UP}_L (f) = \max_{L \in {\cal L}} \mathrm{UP}_L (f),  
$$
and 
${\cal L}$ is a set of all vectors $L\in\rd$ such that $\|L\|=1$, $v:=(L_1^2,\dots,L_d^2)$, and 
$v=B^{-1}E/\|B^{-1} E\|_1,$ where $B \in {\bf A},$ and $\|B^{-1} E\|_1$ is the $l_1$-norm of the vector $B^{-1} E.$   
\end{theo}

In the proof we will show that  ${\cal L}$ is a finite nonempty set, namely $1\leq \#{\cal L}\le d!\sum_{k=1}^d (k!)^{-1}$. So,  Theorem \ref{1depend_L} reduces the extremal problems to calculate  $\mathrm{UP}_L(f)$ for  a finite number of vectors $L$.

If the function $f$ does not meet relations (\ref{symm}) then $\mathrm{UP}_L$ is not a quadratic form anymore and finding its extremal values is a complicated problem allowing  numerical solutions only.     
On the other hand, it turns out that as in the  one-dimensional case the inequalities 
$$
\left(\Delta({\cal A}_L, f)  \Delta({\cal B}_L, f)\right)^{1/2} \geq \frac{1}{2}
C \|f\|_2^2
\ \ \ 
\mbox{ and }\ \ \ 
(2\pi)^{-2}\Delta({\cal A}_L, f) + (2\pi)^2\Delta({\cal B}_L, f) \geq C\|f\|_2^2
$$
are equivalent. 
Indeed, the first inequality implies the second one because of the elementary inequality $2a b \leq a^2+b^2.$ Conversely, substituting the function $c^{d/2}f(c \cdot)$ for $f$ in the second inequality, we get 
$$
(2\pi)^{-2} c^{-2} \Delta({\cal A}_L, f) +(2\pi)^2 c^2 \Delta({\cal B}_L, f)  \geq C \|f\|_2^2
$$
and as 
$c=\left(\Delta({\cal A}_L, f)\right)^{1/4} \left(\Delta({\cal B}_L, f)\right)^{-1/4}$
the last inequality takes the form 
$$2 \left(\Delta({\cal A}_L, f)  \Delta({\cal B}_L, f)\right)^{1/2} \geq~C\|f\|_2^2
$$
that has to be proved.

So the functional $(2\pi)^{-2}\Delta({\cal A}_L, f) + (2\pi)^2 \Delta({\cal B}_L, f)$ can also be used as a measure for a localization of a function.
In contrast to $\mathrm{UP}_L$, the functional $(2\pi)^{-2}\Delta({\cal A}_L, f) + (2\pi)^2\Delta({\cal B}_L, f)$ is always a quadratic form with respect to the coordinates of the vector $L$. 

We still fix a function $f\in {\cal D}({\cal A}_L) \bigcap {\cal D}({\cal B}_L)$, $\|f\|_2=1$ 
and solve the  minimization and maximization problems for the new functional in the next theorem.

\begin{theo}
\label{depend_L}
Let $f\in {\cal D}({\cal A}_L) \bigcap {\cal D}({\cal B}_L)$, $\|f\|_2=1$. The values  
$$
\min_{\|L\|=1} \left((2\pi)^{-2}\Delta({\cal A}_L, f) + (2\pi)^2 \Delta({\cal B}_L, f)\right) \ \ \ 
\mbox{ and }
\ \ \ 
\max_{\|L\|=1} \left((2\pi)^{-2}\Delta({\cal A}_L, f) + (2\pi)^2\Delta({\cal B}_L, f)\right)
$$
are equal to 
the minimal and the maximal eigenvalues of the matrix $M=(M_{k,n})_{k,n=1,\dots,d}$ respectively, where 
\begin{equation}\label{def_M}
\begin{aligned}
M_{k,n}&=  \int_{\rd} x_k x_n |f(x)|^2\,\mathrm{d}x
+ (2\pi)^2\int_{\rd} x_k x_n |\widehat{f}(x)|^2\,\mathrm{d}x\\
-&  \left(\int_{\rd}x_k|f(x)|^2\,\mathrm{d}x \right)^2 \left(\int_{\rd}x_n|f(x)|^2\,\mathrm{d}x \right)^2
- (2\pi)^2\left(\int_{\rd}x_k|\widehat{f}(x)|^2\,\mathrm{d}x \right)^2 \left(\int_{\rd}x_n|\widehat{f}(x)|^2\,\mathrm{d}x \right)^2.
\end{aligned}
\end{equation}
The minimum and the maximum are attained by eigenvectors corresponding to these eigenvalues.  
\end{theo}

\begin{rem} 
\label{indep}
It follows from Theorem \ref{depend_L}  that $(2\pi)^{-2}\Delta({\cal A}_L, f) + (2\pi)^2 \Delta({\cal B}_L, f)$ 
does not depend on $L$ if and only if the matrix $M$ has a unique eigenvalue with  multiplicity $d$, that, since the matrix $M$ is symmetric, is equivalent to $M=\lambda I,$ where $I$ is the identity $d\times d$ matrix and $\lambda$ is the eigenvalue. 
  \end{rem}

  
  \subsection{Time and frequency variances in terms of the Hermite functions}
  
  In \cite{deBruijn}, de Bruijn gives an expression for time and frequency variances in terms of the Fourier-Hermite coefficients. In this subsection we generalize this idea to the multivariate case and variances $\Delta({\cal A}_L, f)$, $\Delta({\cal B}_L, f).$ Without loss of generality, by Lemma \ref{lem:UCRdSTD}, we assume 
  \begin{equation}
       \label{0center}
 \langle {\cal A}_L f,f\rangle=0, \quad  \langle {\cal B}_L f,f\rangle=0. 
  \end{equation}
  So,  
    $$
    \Delta({\cal A}_L, f)=\|{\cal A}_L f\|_2^2 \ \  \mbox{ and  } \ \  \Delta({\cal B}_L, f)=\|{\cal B}_L f\|_2^2.
    $$

The $d$-dimensional Hermite functions
are products of one-dimensional  ones
  $$
  h_{\alpha}(x)
  =h_{\alpha_1}(x_1) h_{\alpha_2}(x_2) \dots h_{\alpha_d}(x_d), \quad \alpha\in\zd_+,
  $$
  where for $k\in\n$, $y\in\mathbb{R}$ we choose the Hermite function in the form (see \cite{Folland1989Book})
  $$
  h_{k}(y) = (-1)^k (2^k k! \sqrt{\pi})^{-\frac 12} {\mathrm e}^{\frac{y^2}{2}} D^k {\mathrm e}^{-y^2}.
  $$

  \begin{theo}
  \label{T:hermite}
  Let $f\in {\cal D}({\cal A}_L) \bigcap {\cal D}({\cal B}_L)$ and 
  let a function
  $f$ be expanded in the Fourier-Hermite series $f=\sum_{\alpha\in\zd_+} c_{\alpha} h_{\alpha}$. Then  
  \begin{equation}
      \label{hermite}
     (2 \pi)^{-2} \|{\cal A}_L f\|_2^2 + (2\pi)^2 \|{\cal B}_L f\|_2^2   =
 \sum_{\alpha \in \zd_+ }
 \left(\left|\sum_{n=1}^d L_n \sqrt{\alpha_n}c_{\alpha_1\dots\alpha_n-1\dots\alpha_d} \right|^2+
 \left|\sum_{n=1}^d L_n \sqrt{\alpha_{n}+1}c_{\alpha_1\dots\alpha_n+1\dots\alpha_d}  \right|^2
 \right),
  \end{equation}
  \begin{equation}
      \label{hermite1}
      \begin{aligned}
      \mathrm{UP}_L (f)=\frac{ \|{\cal A}_L f\|_2^2 
      \|{\cal B}_L f\|_2^2}{\|L\|^4 \|f\|_2^4} = &
      \frac {1}{4 \|L\|^4 \|f\|_2^4}  \sum_{\alpha \in \zd_+ }\left|\sum_{n=1}^d L_n \left(\sqrt{\alpha_n}c_{\alpha_1\dots\alpha_n-1\dots\alpha_d}   +\sqrt{\alpha_{n}+1}c_{\alpha_1\dots\alpha_n+1\dots\alpha_d}  \right)\right|^2 \\
      &
     \times \sum_{\alpha \in \zd_+ }\left|\sum_{n=1}^d L_n \left(\sqrt{\alpha_n}c_{\alpha_1\dots\alpha_n-1\dots\alpha_d}    - \sqrt{\alpha_{n}+1}c_{\alpha_1\dots\alpha_n+1\dots\alpha_d}  \right)\right|^2,
      \end{aligned}
  \end{equation}
  where we put $c_{\alpha_1\dots\alpha_n-1\dots\alpha_d}=0$ for $(\alpha_1,\dots,\alpha_n-1,\dots,\alpha_d)\notin \zd_+$.
  \end{theo}

One can deduce the inequality 
$\displaystyle (2 \pi)^{-2}\|{\cal A}_L f\|_2^2 + (2\pi)^2 \|{\cal B}_L f\|_2^2\geq 
 \|L\|^2 \|f\|_2^2$ and, therefore, the uncertainty principle 
 $
 \displaystyle \mathrm{UP}_L(f)\ge 1/4
 $
 from (\ref{hermite}). 
 Indeed,
 $$
 \sum_{\alpha \in \zd_+ }
 \left(\left|\sum_{n=1}^d L_n \sqrt{\alpha_n}c_{\alpha_1\dots\alpha_n-1\dots\alpha_d} \right|^2+
 \left|\sum_{n=1}^d L_n \sqrt{\alpha_{n}+1}c_{\alpha_1\dots\alpha_n+1\dots\alpha_d}  \right|^2
 \right)
 $$
\begin{eqnarray*}
& \geq  &
 \left|
 \sum_{\alpha \in \zd_+ }
 \left|\sum_{n=1}^d L_n \sqrt{\alpha_n}c_{\alpha_1\dots\alpha_n-1\dots\alpha_d} \right|^2-
 \left|\sum_{n=1}^d L_n \sqrt{\alpha_{n}+1}c_{\alpha_1\dots\alpha_n+1\dots\alpha_d}  \right|^2
 \right| 
\\
& = &
 \left|
 \sum_{\alpha \in \zd_+ } \sum_{n=1}^d L_n^2 \alpha_n
 |c_{\alpha_1\dots\alpha_n-1\dots\alpha_d}|^2 - 
  \sum_{\alpha \in \zd_+ } \sum_{n=1}^d L_n^2 (\alpha_n+1)
 |c_{\alpha_1\dots\alpha_n+1\dots\alpha_d}|^2 \right|
\\
& = &\left|\sum_{\alpha \in \zd_+ } \sum_{n=1}^d L_n^2 |c_{\alpha}|^2\right| = \|L\|^2 \|f\|_2^2.
\end{eqnarray*}
 
 Equality (\ref{hermite}) can also be used to improve the inequality
 $\displaystyle (2\pi)^{-2}\|{\cal A}_L f\|_2^2 + (2\pi)^2 \|{\cal B}_L f\|_2^2\geq 
 \|L\|^2 \|f\|_2^2$ and, in the end, the uncertainty principle
 $
 \displaystyle \mathrm{UP}_L(f)\ge 1/4
 $
 for functions with some kind of symmetry. 
 \begin{lem}
 \label{lem:symm}
 Let  $f\in {\cal D}({\cal A}_L) \bigcap {\cal D}({\cal B}_L)$ and 
 \begin{equation}
     \label{pm_symm}
  f(x_1, \dots, x_k, \dots, x_d) = - f(x_1, \dots, -x_k, \dots, x_d)
 \end{equation}
  for $k=1,\dots,d$, $x\in\mathbb{R}^d$. 
Then  
$$
\displaystyle \mathrm{UP}_L(f)\ge \frac{9}{4}.
$$
 \end{lem}

\section{Proof of statements}


{\bf Proof of Lemma~\ref{lemExactUP}.}
Due to Theorem 3.1 in~\cite{Selig2002UC}  the equality in the uncertainty principle is attained if and only if there exist constants $c_1, c_2, d_1, d_2 \in {\mathbb C}$ with $(|c_1|+|d_1|)(|c_2|+|d_2|)>0$ such that
\begin{equation}
c_1 ({\cal A}_L^* - \overline{a}) f = d_1 ({\cal B}_L - b)f, \quad \text{ and } \quad c_2({\cal A}_L - a)f = d_2 ({\cal B}_L^* - \overline{b}) f
\label{fUCEqual}
\end{equation}
and, either at least one of the constants is zero, or
$\frac{d_1}{c_1} = - \overline{\frac{d_2}{c_2}}.$
Here $a= \frac{\langle{\cal A}_L f,f\rangle}{\|f\|_2^2}$ and
$b= \frac{\langle{\cal B}_L f,f\rangle}{\|f\|_2^2}$.

In our case, since ${\cal A}_L$ and ${\cal B}_L$ are self-adjoint, then $a$ and $b$ are real. Therefore, condition~(\ref{fUCEqual}) is equivalent to 
\begin{equation}
c_1 ({\cal A}_L - a) f = d_1 ({\cal B}_L - b)f, 
\label{fUCEqual2}
\end{equation}
for some  $c_1, d_1 \in {\mathbb C}$, $|c_1|+|d_1|>0,$ and, either at least one of the constants is zero, or $\frac{d_1}{c_1}=-\overline{\frac{d_1}{c_1}}$. 
If $c_1=0$ or $d_1=0$, relation~(\ref{fUCEqual2}) implies that $({\cal B}_L - b)f=0$ or $({\cal A}_L - a) f =0$. In any case $f$ should be zero function.
Assume that $c_1 \neq 0$ and $d_1\neq 0$ and denote $ \mathrm{i}\mu = \frac{d_1}{c_1}$, and $\mu\in{\mathbb R}\setminus\{0\}$. So
$$
({\cal A}_L - a) f = \mathrm{i}\mu ({\cal B}_L - b)f.
$$
Due to Lemma~\ref{lem:UCRdSTD} the value of the uncertainty product does not change, if we replace the function $f$ with the following $g(x) = \mathrm{e}^{2\pi \mathrm{i} \langle \beta,x\rangle} f(x+\alpha),$ where $\langle \alpha, L\rangle=a$, $\langle \beta,L\rangle=b$. But for this function
$$
\frac{\langle{\cal A}_L g,g\rangle}{\|g\|_2^2} = 0, \quad
\frac{\langle{\cal B}_L g,g\rangle}{\|g\|_2^2} = 0.
$$
Thus, without loss of generality, assume that $a=0$ and $b=0$. Now, we need to solve the following equation  $$4 \pi^2 \langle L, x\rangle f(x) = - \mu \frac{\partial f}{\partial L}(x),$$
which is a linear partial\textbf{\textbf{\textbf{}}} differential equation.
Let us rewrite it in another form
$$
\sum_{j=1}^d L_j \frac{\partial f}{\partial x_j}(x) = - \frac{4 \pi^2}{\mu} \langle L, x\rangle f(x).
$$
Using the standard methods of solving such partial differential equation, we combine the additional system of equations
$$
\frac{\mathrm{d} x_1}{L_1} = \frac{\mathrm{d} x_2}{L_2} = \dots = \frac{\mathrm{d} x_d}{L_d} = -\frac{\mu  \mathrm{d} f}{4 \pi^2\langle L, x\rangle f(x)}
$$
and find its $d$ independent first integrals
$$
L_2 x_1 - L_1 x_2 = C_1,\ \  L_3 x_1 - L_1 x_3 = C_2,\ \  \dots, L_d x_1 - L_1 x_d = C_{d-1}
$$
and the last integral can be computed from the following considerations. Since $\|L\|=1$,
$$
-\frac{\mu \mathrm{d} f}{4 \pi^2\langle L, x\rangle f(x)} = 
-\sum_{j=1}^d L_j^2 \frac{\mu \mathrm{d} f}{4 \pi^2\langle L, x\rangle f(x)} = \sum_{j=1}^d L_j \mathrm{d} x_j = \mathrm{d} \langle L,x\rangle. 
$$
So,
$$
\frac{ \mathrm{d} f}{f(x)} = -\frac{4 \pi^2 }{\mu} \mathrm{d} \frac{\langle L,x\rangle^2}{2}. 
$$
Therefore, the last first integral is given by
$$
f(x) = C_d \  \mathrm{e}^{-\frac{2 \pi^2 }{\mu}\langle L,x\rangle^2 }.
$$
Since the function $f$ appears only in one first integral, then the general solution can be written as
$$
f(x) = \mathrm{e}^{-\frac{2 \pi^2 }{\mu}\langle L,x\rangle^2 } \Phi (L_2 x_1 - L_1 x_2,L_3 x_1 - L_1 x_3,\dots,L_d x_1 - L_1 x_d),
$$
where $\Phi$ is an arbitrary continuously differentiable function (such that $\mathrm{UP}_L(f)$ makes sense).
For this class of functions $\mathrm{UP}_L(f) = \frac{1}{4}.$ \hfill $\Box$

In order to prove Theorems~\ref{theoL} and~\ref{theoGG} we need some additional statements and notations.
For an admissible function $f$ and a parameter $\lambda \in {\mathbb R}$, we denote $f_\lambda(x) : = \sqrt{\lambda^d} f(\lambda x).$ The function $f_\lambda$ is also admissible. Although, $f_\lambda$ is not periodic, we will use notations $\|f_\lambda\|_{\td}^2 := \int_{\td} |f_\lambda|^2 $ and $\langle f_\lambda,g\rangle_{\td} = \int_{\td} f_\lambda \overline{g}$, assuming that $\td = [-1/2,1/2)^d$, where $g$ is in $L_2(\td)$ or also is an admissible function.

Now, we rewrite  $\langle {\cal A}^{\td}_L f,f\rangle_{\td}$ for an admissible function $f$.
Define two functionals  
$$
K_L(f) = \frac{1}{2} \int_{\td} \left|\mathrm{e}^{2\pi \mathrm{i} \langle L, x\rangle}-1\right|^2 |f(x)|^2 \mathrm{d}x = 2 \int_{\td} \sin^2 \frac{2\pi \langle L, x \rangle}{2} |f(x)|^2 \mathrm{d}x,
$$
$$
M_L(f) = \frac{1}{2} \int_{\td} (\mathrm{e}^{2\pi \mathrm{i} \langle L, x\rangle}-1) (\mathrm{e}^{-2\pi \mathrm{i} \langle L, x\rangle} + 1) |f(x)|^2 \mathrm{d}x = \mathrm{i} \int_{\td} \sin (2\pi \langle L, x \rangle) |f(x)|^2 \mathrm{d}x.
$$

From 
$$
2 \|f\|^2_{\td} - 2 \ {\cal R}e(\langle {\cal A}^{\td}_L f,f\rangle_{\td}) = \int_{\td} (2 - \mathrm{e}^{2\pi \mathrm{i} \langle L, x\rangle} - \mathrm{e}^{-2\pi \mathrm{i} \langle L, x\rangle}) |f(x)|^2 \mathrm{d}x = 2 K_L(f)
$$
it follows that ${\cal R}e(\langle {\cal A}^{\td}_L f,f\rangle_{\td}) = \|f\|^2_{\td} - K_L(f)$.
Also,
$$
{\cal I}m(\langle {\cal A}^{\td}_L f,f\rangle_{\td})= \frac{1}{2 \mathrm{i}} \int_{\td} (\mathrm{e}^{2\pi \mathrm{i} \langle L, x\rangle}-1) (\mathrm{e}^{-2\pi \mathrm{i} \langle L, x\rangle} + 1) |f(x)|^2 \mathrm{d}x = - \mathrm{i} M_L(f),
$$
since $\mathrm{e}^{2\pi \mathrm{i} \langle L, x\rangle}-\mathrm{e}^{-2\pi \mathrm{i} \langle L, x\rangle} = (\mathrm{e}^{2\pi \mathrm{i} \langle L, x\rangle}-1)(\mathrm{e}^{-2\pi \mathrm{i} \langle L, x\rangle}+1)$.
Thus,
$$
|\langle {\cal A}^{\td}_L f,f\rangle_{\td}|^2 = (\|f\|_{\td}^2 - K_L(f))^2 - M_L^2(f).
$$
The directional angular variance can be written as follows
\begin{equation}
    \text{var}_L^{\mathrm A} (f) = \frac{\|f\|_{\td}^4 - (\|f\|_{\td}^2 - K_L(f))^2 + M_L^2(f)}{(\|f\|_{\td}^2 - K_L(f))^2 - M_L^2(f)} = \frac{2 \|f\|_{\td}^2 K_L(f) - K_L^2(f) + M_L^2(f)}{(\|f\|_{\td}^2 - K_L(f))^2 - M_L^2(f)}.
\label{eq:DirAngVar}
\end{equation}

\begin{lem}
    Let $f$ be an admissible function and $\lambda>0,$
    $L\in\zd.$ Then
    $$
    \lim_{\lambda\to\infty} \|f_\lambda\|^2_{\td} = \|f\|^2_2,\quad\quad
    \lim_{\lambda\to\infty} \frac{1}{\lambda^2} \left\|{\cal B}_Lf_\lambda\right\|_{\td}^2 = \left\|{\cal B}_Lf\right\|_2^2,
    $$
    $$
    \lim_{\lambda\to\infty} \frac{1}{\lambda} \left\langle {\cal B}_Lf_\lambda, f_\lambda \right\rangle_{\td} = \left\langle {\cal B}_Lf, f \right\rangle.
    $$
Additionally,
$$
\lim_{\lambda\to\infty} 2 \lambda^2 K_L(f_\lambda) = \|{\cal A}_Lf\|^2_2, \qquad
\lim_{\lambda\to\infty} \lambda M_L(f_\lambda) = i \langle {\cal A}_Lf, f \rangle.
$$
\label{lem:UCTDtoRD}
\end{lem}
The proof can be given by straightforward computations following the proof of the analogous results in~\cite{PrestinQRS2003UC}.

Now, we study the behavior of the periodized version of a scaled admissible function, i.e.
$
f_\lambda^{\mathrm{per}} (x) 
=  \sum_{k\in\zd} f_\lambda(x+k). 
$

\begin{lem}
    Let $f$ be an admissible function, $L\in\zd$ and $\lambda>0.$ Then $f_\lambda^{\mathrm{per}}$ and $\frac{\partial f_\lambda^{\mathrm{per}}}{\partial L}$ are continuous functions in $L_2(\td).$
\end{lem}

{\bf Proof.} For $x\in\td$, we get the following estimate for a big enough $N\in{\mathbb N}$, using the admissibility of $f$,
$$
\left| \sum_{k\in\zd} f_\lambda(x+k) - \sum_{\|k\|<N} f_\lambda(x+k)\right| 
\le \sum_{\|k\|\ge N} \sqrt{\lambda^d} \left|f(\lambda(x+k))\right|
$$
$$
 \le 
  \sum_{\|k\|\ge N} \frac{C_1 \sqrt{\lambda^d}}{ \|\lambda(x+k)\|^\gamma} \le  \sum_{\|k\|\ge N} \frac{C_1 \lambda^{d/2 - \gamma}}{(\|k\|-\|x\|)^\gamma}
\le \sum_{\|k\|\ge N} \frac{C_1 \lambda^{d/2 - \gamma}}{(\|k\|-\sqrt{d}/2)^\gamma} \to 0, \quad N\to\infty.
$$
The last expression is independent of $x$. Since $f_\lambda$ is continuous and the convergence of the series $\sum_{k\in\zd} f_\lambda(x+k)$ is uniform, $f_\lambda^{\mathrm{per}}$ is also continuous.
The same estimate is valid for the continuous function $\frac{\partial f_\lambda}{\partial L}$, and therefore, $\frac{\partial f_\lambda^{\mathrm{per}}}{\partial L}$ is a continuous function. \hfill $\Box$

Now, we study the limit behavior of the directional angular and frequency variances of $f_\lambda^{\mathrm{per}}.$ 
\begin{lem}
    Let $f$ be admissible, $L\in\zd$ and $\lambda>0.$ Then
    $$
    \lim_{\lambda\to\infty} \frac{1}{\lambda^2} \ \text{\rm var}_L^{\mathrm F}(f_\lambda^{\mathrm{per}}) = \frac{\Delta({\cal B}_L,f)}{\|f\|_2^2}.
    $$
    \label{lem:UCTDtoRDfreq}
\end{lem}

{\bf Proof.} Using the admissibility of $f$, for $k\neq \nul$ we get
$$
\| f_\lambda(\cdot + k)\|^2_{\td} = \lambda^d
\int_{\td} |f(\lambda(x + k))|^2 \mathrm{d}x 
\le \int_{\td} \frac{C_1^2 \lambda^d}{\|\lambda(x + k)\|^{2\gamma}} \mathrm{d}x \le \int_{\td} \frac{C_1^2 \lambda^{d-2\gamma}}{\|x + k\|^{2\gamma}} \mathrm{d}x.
$$
For big enough $k$ ($\|k\|> \sqrt{d}/2$), 
$$
\frac{1}{\|x + k\|^{2\gamma}} \le \frac{1}{(\| k\| - \| x\|)^{2\gamma}} \le \frac{1}{(\| k\| - \sqrt{d}/2)^{2\gamma}}.
$$
Therefore, we can state that $S(\gamma):=\sum\limits_{k\neq \nul}  \left(\int_{\td} \|x + k\|^{-2\gamma} \mathrm{d}x \right)^{1/2} < \infty$ and
$$
\sum_{k\neq \nul} \|f_\lambda(\cdot + k)\|_{\td} \le C_1 \lambda^{d/2-\gamma}  \sum_{k\neq \nul}  \left(\int_{\td} \frac{\mathrm{d}x}{\|x + k\|^{2\gamma}} \right)^{1/2}= C_1 \lambda^{d/2-\gamma}  S(\gamma) \to 0, \quad \lambda\to\infty.
$$
Analogously, we can estimate the derivatives $D^{e_j} f_\lambda$, $j=1,\dots,d$. Namely,
$$
\|D^{e_j} f_\lambda (\cdot + k)\|^2_{\td} = \int_{\td} |D^{e_j} f_\lambda (x + k)|^2 \mathrm{d}x = \lambda^{d+2} \int_{\td} |D^{e_j} f(\lambda(x + k))|^2 \mathrm{d}x \le
\int_{\td} \frac{C_2^2 \lambda^{d+2-2\beta}}{\|x + k\|^{2\beta}} \mathrm{d}x.
$$
Therefore,
$$
\left \|\frac{\partial f_\lambda}{\partial L} (\cdot + k)\right\|^2_{\td} \le \sum_{j=1}^d L^2_j \|D^{e_j} f_\lambda (\cdot + k)\|^2_{\td}  \le \|L\|^2
\int_{\td} \frac{C_2^2 \lambda^{d+2-2\beta}}{\|x + k\|^{2\beta}} \mathrm{d}x.
$$
Hence,
$$
\frac{1}{\lambda} \sum_{k\neq \nul} \left \|\frac{\partial f_\lambda}{\partial L} (\cdot + k)\right\|_{\td}  \le C_2 \|L\| \lambda^{d/2-\beta}  \sum_{k\neq \nul}  \left(\int_{\td} \frac{\mathrm{d}x}{\|x + k\|^{2\beta}} \right)^{1/2} = C_2 \|L\| \lambda^{d/2-\beta}  S(\beta) \to 0, 
$$
as $\lambda\to\infty.$ Now, consider, 
$$
\left|\ \|f_\lambda^{\mathrm{per}}\|_{\td} - \|f\|_2 \right| \le \left| \|f_\lambda^{\mathrm{per}}\|_{\td} - \|f_\lambda\|_{\td} \right| + \left| \|f_\lambda\|_{\td} -  \|f\|_2 \right|.
$$
The first term can be estimated as follows
\begin{eqnarray*}
\left| \ \|f_\lambda^{\mathrm{per}}\|_{\td} - \|f_\lambda\|_{\td} \right| &=& \left| \left\|f_\lambda + \sum_{k\neq \nul} f_\lambda(\cdot + k)\right\|_{\td} - \|f_\lambda\|_{\td} \right|
\\
 &\le & \left\|\sum_{k\neq \nul} f_\lambda(\cdot + k)\right\|_{\td} \le \sum_{k\neq \nul} \left\| f_\lambda(\cdot + k)\right\|_{\td} \to 0, \quad  \lambda\to\infty.
\end{eqnarray*}
The second term tends to zero as $\lambda\to\infty$ by Lemma~\ref{lem:UCTDtoRD}.
Thus, we get
$$
\lim_{\lambda\to\infty} \|f_\lambda^{\mathrm{per}}\|^2_{\td} = \|f\|^2_2.
$$
Analogously, it can be stated that
$$
\lim_{\lambda\to\infty} \frac{1}{\lambda^2} \left\|\frac{\partial f_\lambda^{\mathrm{per}}}{\partial L}\right\|^2_{\td} = \left\|\frac{\partial f}{\partial L}\right\|^2_2.
$$
Furthermore,
\begin{eqnarray*}
\frac{1}{\lambda} \left\langle \frac{\partial f_\lambda^{\mathrm{per}}}{\partial L}, f_\lambda^{\mathrm{per}} \right\rangle_{\td} &= &
\frac{1}{\lambda}  \left\langle \sum_{k\in\zd}\frac{\partial f_\lambda}{\partial L} (\cdot + k), \sum_{l\in\zd} f_\lambda(\cdot + l) \right\rangle_{\td} 
\\
&= & \frac{1}{\lambda} \left\langle \frac{\partial f_\lambda}{\partial L}, f_\lambda \right\rangle_{\td} + \frac{1}{\lambda} \sum_{l\in\zd, l\neq \nul} \left\langle \frac{\partial f_\lambda}{\partial L}, f_\lambda(\cdot + l) \right\rangle_{\td} 
\\
& &+\frac{1}{\lambda} \sum_{k\in\zd, k\neq \nul} \sum_{l\in\zd} \left\langle \frac{\partial f_\lambda}{\partial L}(\cdot+k), f_\lambda(\cdot + l) \right\rangle_{\td}.
\end{eqnarray*}
With the Cauchy-Bunyakovsky-Schwarz inequality and above considerations, we estimate the last two terms as
$$
\left| \frac{1}{\lambda} \sum_{l\in\zd, l\neq \nul} \left\langle \frac{\partial f_\lambda}{\partial L}, f_\lambda(\cdot + l) \right\rangle_{\td} \right| \le
\frac{1}{\lambda} \left\| \frac{\partial f_\lambda}{\partial L} \right\|_{\td} \sum_{l\neq \nul} \|f_\lambda(\cdot + l)\|_{\td} \to 0, \quad \lambda \to\infty,
$$
since Lemma~\ref{lem:UCTDtoRD} states that $\frac{1}{\lambda} \left\| \frac{\partial f_\lambda}{\partial L} \right\|_{\td} \to \left\| \frac{\partial f_\lambda}{\partial L} \right\|_2 < \infty$, and
$$
\left| \frac{1}{\lambda}  \sum_{k\in\zd, k\neq \nul} \sum_{l\in\zd} \left\langle \frac{\partial f_\lambda}{\partial L}(\cdot+k), f_\lambda(\cdot + l) \right\rangle_{\td}\right|  \hspace{8cm}
$$
$$
 \le \sum_{k\in\zd, k\neq \nul} \frac{1}{\lambda} \left\| \frac{\partial f_\lambda}{\partial L}(\cdot+k) \right\|_{\td}
\sum_{l\in\zd} \|f_\lambda(\cdot + l)\|_{\td} \to 0, \quad \lambda\to 0,
$$
since 
by 
$\|f_\lambda\|_{\td} \to \|f\|_2<\infty.$
Thus,
$$
\lim_{\lambda\to\infty}\frac{1}{\lambda} \left\langle \frac{\partial f_\lambda^{\mathrm{per}}}{\partial L}, f_\lambda^{\mathrm{per}} \right\rangle_{\td} = \left\langle \frac{\partial f}{\partial L}, f \right\rangle.
$$
Combining all the limits together and noting that
$$
\lim_{\lambda\to\infty} \frac{1}{\lambda^2} \left\| {\cal B}_L^{\td} f_\lambda^{\mathrm{per}}\right\|^2_{\td} = \left\|{\cal B}_L f\right\|^2_2, \qquad
\lim_{\lambda\to\infty}\frac{1}{\lambda} \left\langle {\cal B}_L^{\td} f_\lambda^{\mathrm{per}}, f_\lambda^{\mathrm{per}} \right\rangle_{\td} = \left\langle {\cal B}_L f , f \right\rangle
$$
we get
$$
\lim_{\lambda\to\infty} \frac{1}{\lambda^2} \text{var}{\mathrm F}_L(f_\lambda^{\mathrm{per}}) = 
\lim_{\lambda\to\infty} \left( \frac{\frac{1}{\lambda^2} \left\|{\cal B}_L^{\td} f_\lambda^{\mathrm{per}}\right\|^2_{\td}}{\|f_\lambda^{\mathrm{per}}\|_{\td}^2}  - \frac{\frac{1}{\lambda^2} \left|\left\langle {\cal B}_L^{\td} f_\lambda^{\mathrm{per}}, f_\lambda^{\mathrm{per}} \right\rangle_{\td}\right|^2}{\|f_\lambda^{\mathrm{per}}\|_{\td}^4}\right) = \frac{\Delta({\cal B}_L,f)}{\|f\|_2^2}.
$$
\hfill $\Box$

Now we consider the directional angular variance.
\begin{lem}
    Let $f$ be admissible, $L\in\zd$ and $\lambda>0.$ Then there exists a $\lambda_1>0,$ such that $\text{\rm var}_L{\mathrm A}(f_\lambda^{\mathrm{per}})$ is finite for all $\lambda>\lambda_1$ and
    $$
    \lim_{\lambda\to\infty} \lambda^2 \text{\rm var}_L{\mathrm A}(f_\lambda^{\mathrm{per}}) = \frac{\Delta({\cal A}_L,f)}{\|f\|_2^2}. 
    $$
    \label{lem:VarATdRd}
\end{lem}

{\bf Proof.} We will use the representation of the directional angular variance~(\ref{eq:DirAngVar}).
Let us consider $K_L(f_\lambda^{\mathrm{per}})$ first. Since $\sin^2 \frac{2\pi \langle L, x\rangle}{2}\ge 0$ and $\sin^2 \frac{2\pi \langle L, x\rangle}{2} = 0$ on a set of measure zero,  $\sqrt{K_L(\cdot)}$ is actually a weighted $L_2$ norm. This allows to proceed as follows.
Note that 
\begin{eqnarray*}
K_L(f_\lambda(\cdot+k)) &=& 2 \int_{\td} \sin^2 \frac{2\pi \langle L, x \rangle}{2} |f_\lambda(x+k)|^2 \mathrm{d}x 
\\
&\le &  2 \int_{\td}  |f_\lambda(x+k)|^2 \mathrm{d}x = 2 \|f_\lambda(\cdot + k)\|^2_{\td}.
\end{eqnarray*}
By the admissibility of $f_\lambda$ and the estimates in Lemma~\ref{lem:UCTDtoRDfreq}, we obtain
$$
\sqrt{2} \lambda \sum_{k\neq \nul} \sqrt{K_L(f_\lambda(\cdot+k))} \le 2 \lambda \sum_{k\neq \nul} \|f_\lambda(\cdot + k)\|_{\td} \le  2  C \lambda^{d/2 +1 -\gamma}  S(\gamma) \to 0, \quad \lambda\to\infty,
$$
since $\gamma> d/2+1$.
Using the triangle inequality for the weighted norm, we get
$$
\sqrt{2} \lambda \left| \sqrt{K_L(f_\lambda^{\mathrm{per}})} - \sqrt{K_L(f_\lambda)}\right| \le \sqrt{2} \lambda  \sum_{k\neq \nul} \sqrt{K_L(f_\lambda(\cdot+k))}.
$$
Now, it can be stated that
\begin{equation}
\lim_{\lambda\to\infty} 2 \lambda^2 K_L(f_\lambda^{\mathrm{per}}) = \|{\cal A}_L f\|_2^2.
\label{fKL}
\end{equation}
Indeed, since
$$
\left|\sqrt{2} \lambda \sqrt{K_L(f_\lambda^{\mathrm{per}})} - \|{\cal A}_L f\|_2  \right| \le 
\left|\sqrt{2} \lambda \sqrt{K_L(f_\lambda^{\mathrm{per}})} -\sqrt{2} \lambda \sqrt{K_L(f_\lambda)}  \right|  + \left| \sqrt{2} \lambda \sqrt{K_L(f_\lambda)} -  \|{\cal A}_L f\|_2  \right|,
$$
where the first term goes to zero as $\lambda\to\infty$ by the above inequality, the second one by Lemma~\ref{lem:UCTDtoRD}.

Now, we establish that 
$\lim_{\lambda\to\infty} \lambda M_L(f^{\mathrm{per}}_\lambda) = \mathrm{i} \langle {\cal A}_Lf, f \rangle.
$
Again, we start from the following estimate 
$$
\left|\lambda M_L(f^{\mathrm{per}}_\lambda) - \mathrm{i} \langle {\cal A}_Lf, f \rangle \right| \le \left|  \lambda M_L(f^{\mathrm{per}}_\lambda) - \lambda M_L(f_\lambda)\right| + \left|  \lambda M_L(f_\lambda) -  \mathrm{i} \langle {\cal A}_Lf, f \rangle \right|.
$$
By Lemma~\ref{lem:UCTDtoRD} the second term tends to zero since
$\lim_{\lambda\to\infty} \lambda M_L(f_\lambda) = \mathrm{i} \langle {\cal A}_Lf, f \rangle.$
Thus, it is sufficient to prove, that
$$
\lim_{\lambda\to\infty}  \lambda M_L(f^{\mathrm{per}}_\lambda) = \lambda M_L(f_\lambda).
$$
Recall that for some admissible $g$ $$
M_L(g) =  \mathrm{i} \int_{\td} \sin (2\pi \langle L, x \rangle) |g(x)|^2 \mathrm{d}x.
$$
So, $\sqrt{M_L(\cdot)}$ is not a weighted norm. But it is possible to split the area of integration as follows. Let
$$
P^+ = \{ x\in\td, \ \sin  \langle 2\pi L, x \rangle > 0\}, \quad P^- = \{ x\in\td, \ \sin  \langle 2\pi L, x \rangle < 0\}
$$
and
 $$
M_L^+(g) =   \int_{P^+} \sin (2\pi \langle L, x \rangle) |g(x)|^2  \mathrm{d}x, \quad 
M_L^-(g) =   \int_{P^-} (-\sin 2\pi \langle L, x \rangle) |g(x)|^2 \mathrm{d}x.
$$
Then, $
M_L(g) = \mathrm{i} (M_L^+(g) - M_L^-(g)).
$
Thus, $\sqrt{M_L^+(\cdot)}$ is a weighted norm on  $L_2 (P^+)$, $\sqrt{M_L^-(\cdot)}$ is a weighted norm on  $L_2 (P^-)$.
Therefore, it is sufficient to prove that
$$
\lim_{\lambda\to\infty}  \lambda M_L^+(f^{\mathrm{per}}_\lambda) = \lambda M_L^+(f_\lambda), \quad \lim_{\lambda\to\infty}  \lambda M_L^-(f^{\mathrm{per}}_\lambda) = \lambda M_L^-(f_\lambda).
$$
Hence,
$$
\left|\sqrt{\lambda} \sqrt{M_L^+(f^{\mathrm{per}}_\lambda)} - \sqrt{\lambda} \sqrt{M_L^+(f_\lambda)} \right| \le \left|\sqrt{\lambda} \sqrt{M_L^+\left(\sum_{k\neq\nul} f_\lambda(\cdot+k)\right)}  \right| \le \sqrt{\lambda} \sum_{k\neq\nul} \sqrt{M_L^+\left(f_\lambda(\cdot+k)\right)}.
$$
Consider the following estimates
$$
M_L^+\left( f_\lambda(\cdot+k)\right) = \int_{P^+} \sin (2\pi \langle L, x \rangle) |f_\lambda(x+k)|^2 \mathrm{d}x \le  \|f_\lambda(\cdot + k)\|^2_{\td},
$$
$$
\sqrt{\lambda} \sum_{k\neq\nul} \sqrt{M_L^+\left(f_\lambda(\cdot+k)\right)} \le \sqrt{\lambda} \sum_{k\neq\nul} \|f_\lambda(\cdot + k)\|_{\td} \le  2  C \lambda^{d/2 +1/2 -\gamma}  S(\gamma) \to 0, \quad \lambda\to\infty
$$
and their counterparts for $M_L^-$.
Thus,
\begin{equation}
  \lim_{\lambda\to\infty} \lambda M_L(f^{\mathrm{per}}_\lambda) = \mathrm{i} \langle {\cal A}_Lf, f \rangle.
  \label{fML}
\end{equation}
Recall that
$$
|\langle {\cal A}^{\td}_L f,f\rangle_{\td}|^2 = (\|f\|_2^2 - K_L(f))^2 - M_L^2(f).
$$
Therefore,
$$
\lim_{\lambda\to\infty}  |\langle {\cal A}^{\td}_L f^{\mathrm{per}}_\lambda,f^{\mathrm{per}}_\lambda\rangle_{\td}|^2 = \lim_{\lambda\to\infty}  (\|f^{\mathrm{per}}_\lambda\|_2^2 - K_L(f^{\mathrm{per}}_\lambda))^2 - M_L^2(f^{\mathrm{per}}_\lambda) = \|f\|_2^4 >0,
$$
since $M_L(f^{\mathrm{per}}_\lambda)$ and $K_L(f^{\mathrm{per}}_\lambda)$ should tend to zero as $\lambda\to\infty$ by (\ref{fKL}) and (\ref{fML}). Also, there exists  big enough $\lambda_1$, such that for any $\lambda> \lambda_1$, $(\|f^{\mathrm{per}}_\lambda\|_2^2 - K_L(f^{\mathrm{per}}_\lambda))^2 - M_L^2(f^{\mathrm{per}}_\lambda)>0$.
Thus,
$$
\lim_{\lambda\to\infty}  \lambda^2 \text{var}_L{\mathrm A} (f^{\mathrm{per}}_\lambda) = \lim_{\lambda\to\infty}  \frac{2 \lambda^2 \|f^{\mathrm{per}}_\lambda\|_{\td}^2 K_L(f^{\mathrm{per}}_\lambda) - \lambda^2 K_L^2(f^{\mathrm{per}}_\lambda) + \lambda^2 M_L^2(f^{\mathrm{per}}_\lambda)}{(\|f^{\mathrm{per}}_\lambda\|_{\td}^2 - K_L(f^{\mathrm{per}}_\lambda))^2 - M_L^2(f^{\mathrm{per}}_\lambda)}  $$
$$= \frac{\|f\|_2^2 \|{\cal A}_L f\|^2_2 - \langle {\cal A}_Lf, f \rangle^2 }{\|f\|_2^4} = \frac{\Delta({\cal A}_L,f)}{\|f\|_2^2}.
$$
The last equality is valid since ${\cal A}_L$ is self-adjoint and therefore, $\langle {\cal A}_Lf, f \rangle$ is real.
\hfill $\Box$


{\bf Proof of Theorem~\ref{theoL}.} To prove the connection between $\mathrm{UP}_{L}^{\td}$ and $\mathrm{UP}_{L}$, which are defined in~(\ref{fUPTD}) and~(\ref{fUPRD}), we  apply Lemmas~\ref{lem:UCTDtoRDfreq},~\ref{lem:VarATdRd} to~(\ref{fUPTD}) and get that $\mathrm{UP}_L^{\td} (f_\lambda^{\mathrm{per}}) \to \mathrm{UP}_L(f)$ as $\lambda \to +\infty.$
\hfill $\Box$
    
{\bf  Proof of Theorem~\ref{theoGG}.}
  To prove the connection between $\mathrm{UP_{GG}^{\td}}$ and $\mathrm{UP_{GG}}$
  which are defined in~(\ref{fUPGG}) and~(\ref{fUPGGTD})
  we can use Lemmas~\ref{lem:UCTDtoRDfreq} and~\ref{lem:VarATdRd} with $L=e_j.$
Namely, it is straightforward to see by Lemma~\ref{lem:UCTDtoRDfreq} that
    $$
    \lim_{\lambda\to\infty} \frac{1}{\lambda^2} \ \text{\rm var}_{\mathrm{GG}}^{\mathrm F}(f_\lambda^{\mathrm{per}}) = \sum\limits_{j=1}^d \frac{\Delta({\cal B}_j,f)}{\|f\|_2^2}.
    $$
Also by the proof of Lemma~\ref{lem:VarATdRd} it can be shown that
\begin{eqnarray*}
    \lim_{\lambda\to\infty} \lambda^2 \text{\rm var}_{\mathrm{GG}}^{\mathrm A}(f_\lambda^{\mathrm{per}}) &= &
     \lim_{\lambda\to\infty} \frac{ \sum\limits_{j=1}^d \left(\lambda^2\|f_\lambda^{\mathrm{per}}\|^4_{\td} - \lambda^2|\langle{\cal A}_j^{\td} f_\lambda^{\mathrm{per}},f_\lambda^{\mathrm{per}}\rangle_{\td}|^2\right)} {\left( \sum\limits_{j=1}^d |\langle{\cal A}_j^{\td} f_\lambda^{\mathrm{per}},f_\lambda^{\mathrm{per}}\rangle_{\td}|\right)^2} 
     \\
     &=& \frac{\sum\limits_{j=1}^d \|f\|_2^2 \|{\cal A}_j f\|^2_2 - \langle {\cal A}_jf, f \rangle^2 }{d^2\|f\|_2^4} = \frac{\sum\limits_{j=1}^d \Delta({\cal A}_j,f)}{d^2\|f\|_2^2}. 
\end{eqnarray*}
    Therefore, $\mathrm{UP_{GG}}^{\td} (f_\lambda^{\mathrm{per}}) \to \mathrm{UP_{GG}}(f)$ as $\lambda \to +\infty.$
\hfill $\Box$

{\bf Proof of Theorem \ref{1depend_L}.} 
Writing the time variance in detail we obtain
\begin{equation}
\begin{aligned}
\Delta({\cal A}_L, f) &= \|{\cal A}_L f\|_2^2 - \left|\langle{\cal A}_L f,f\rangle\right|^2 \\
&= 4\pi^2 \int\limits_{\rd} \left|\sum_{j=1}^{d}L_j x_j\right|^2  |f(x)|^2\, \mathrm{d}x - 4\pi^2 \left(\,\,\int\limits_{\rd}\sum_{j=1}^{d}L_j x_j\cdot |f(x)|^2\, \mathrm{d}x\right)^2.
\end{aligned}
\label{varX}
\end{equation}

Similarly, using the property of the Fourier transform we get for the frequency variance
\begin{equation}
\begin{aligned}
\Delta({\cal B}_L, f) &= \|{\cal B}_L f\|_2^2 - \left|\langle{\cal B}_L f,f\rangle\right|^2 \\
&=\int\limits_{\rd}  \left|\frac{ \mathrm{i} }{2\pi} \frac{\partial f}{\partial L} (x)\right|^2 \, \mathrm{d}x - \left(\int\limits_{\rd}\frac{\mathrm{i}}{2\pi}\frac{\partial f}{\partial L}\cdot \overline{f(x)}\, \mathrm{d}x\right)^2\\
&=\int\limits_{\rd} \left|\sum_{j=1}^{d}L_j x_j\right|^2  |\widehat{f}(x)|^2\, \mathrm{d}x - \left(\,\,\int\limits_{\rd}\sum_{j=1}^{d}L_j x_j\cdot |\widehat{f}(x)|^2\, \mathrm{d}x\right)^2.
\end{aligned}
\label{varF}
\end{equation}
Since $|f|$ and $|\widehat{f}|$ are even with respect to each variable (see~(\ref{symm})), then for all $k,j=1,\dots,d$, $k\neq j,$
$$
0 =  \int_{\rd} x_k |f(x)|\, \mathrm{d}x = \int_{\rd}x_k x_j |f(x)|\, \mathrm{d}x  = \int_{\rd} x_k |\widehat{f}(x)|\, \mathrm{d}x = \int_{\rd}x_k x_j |\widehat{f}(x)|\, \mathrm{d}x.
$$
So, $\mathrm{UP}_L$ takes the form
$$
\mathrm{UP}_L  =  \sum_{k=1}^d L_k^2 M_k \sum_{k=1}^d L_k^2 \widehat{M}_k = v^{\mathrm T} A v,
$$
where $A$ is a $d\times d$ matrix whose elements are $(M_k\widehat{M_j}+M_j\widehat{M_k})/2$, $j, k =1,\dots,d$, and 
$v:=(L_1^2,\dots,L_d^2).$

Therefore, to find 
${\rm min}_{\|L\|=1} \mathrm{UP}_L (f)$ and ${\rm max}_{\|L\|=1} \mathrm{UP}_L (f)$
we derive at the following extremal problem with respect to the vector $v$  for the quadratic form  
$v^{\mathrm{T}}Av$ 
$$
\left\{
\begin{array}{l}
v^{\mathrm{T}}Av \to {\rm extr},\\
v_1+\dots+v_d =1, \ \ v_j\ge 0, \ \  j=1,\dots,d. 
\end{array}
\right.
$$
Since the restriction set  $ V:=\{v\in\rd \, ;\, v_1+\dots,+v_d =1, v_j\ge 0,  j=1,\dots,d\}$ is compact, it follows that 
the solution for the extremal problem exists. It remains to follow the well-known classical scheme for a solution of such problems. According to this scheme, extremal points lay on the boundary of the restriction set $V$ or they  are contained among the solutions of the systems of equations 
$$
\frac{\partial}{\partial v_j}\left(v A v^{\mathrm{T}} - \lambda (v_1+\dots +v_d)\right) = 0,  \ \ j=1\dots,d, \ \ \lambda \in \mathbb{R}.
$$
The last system is rewritten in the form 
$
2 A v = \lambda E,
$ 
where $E=(1,\dots,1)\in\rd.$ Thus, 
$v = \lambda/2 A^{-1} E,$ and the Lagrange parameter $\lambda$ is chosen to meet the condition $v_1+ \dots +v_d=1$.  
If extremal points lay on the boundary of the set $V$ then we come to the analogous system of equations, however the matrix $A$ is replaced by the matrix $A_{j_1,\dots,j_q}.$ \hfill $\Box$

{\bf Proof of Theorem \ref{depend_L}.} Starting with formulas (\ref{varX}) and (\ref{varF}) we obtain 
the following quadratic form
$$
(2\pi)^{-2}\Delta({\cal A}_L, f) + (2\pi)^2\Delta({\cal B}_L, f) = L^{\mathrm T} M L,
$$
where the matrix $M$ is defined by \eqref{def_M}.
So, the statement of the theorem is a  well-known fact of linear algebra. \hfill $\Box$

 {\bf Proof of Theorem \ref{T:hermite}.}
  Using the recurrent formula (see \cite{Folland1989Book}) 
  $$
  \sqrt{2} x_n h_{\alpha_n}(x_n) =\sqrt{\alpha_n+1} h_{\alpha_{n+1}}(x_n)+\sqrt{\alpha_n} h_{\alpha_{n-1}}(x_{n}), \quad
 n=1,\dots,d,
 $$
  we obtain
\begin{eqnarray*}
  (2\pi)^{-1}{\cal A}_L f &=&  \sum_{\alpha \in \zd_+ }c_{\alpha}
  \sum_{n=1}^d L_n x_n h_{\alpha_1}(x_1) h_{\alpha_2}(x_2) \dots h_{\alpha_d}(x_d)
\\
  &= &
  \frac{1}{\sqrt{2}}
  \sum_{\alpha \in \zd_+ }c_{\alpha}
  \sum_{n=1}^d L_n \left(\sqrt{\alpha_n+1} h_{\alpha_{n+1}}(x_n)+\sqrt{\alpha_n} h_{\alpha_{n-1}}(x_{n}) \right)
\\
 & & \hspace{4cm} \times
  h_{\alpha_1}(x_1)\dots h_{\alpha_{n-1}}(x_{n-1})
  h_{\alpha_{n+1}}(x_{n+1})
  \dots h_{\alpha_d}(x_d)
\\
  &=&
   \frac{1}{ \sqrt{2}} 
  \sum_{\alpha \in \zd_+ } h_{\alpha}(x) 
  \sum_{n=1}^d L_n \left(\sqrt{\alpha_n}c_{\alpha_1\dots\alpha_n-1\dots\alpha_d}   +\sqrt{\alpha_{n}+1}c_{\alpha_1\dots\alpha_n+1\dots\alpha_d}  \right). 
\end{eqnarray*}
 Here we set $h_{\alpha}(x)=0$ and $c_{\alpha}=0$ for $\alpha \notin \zd_+.$
 Thus, due to orthonormality of the Hermite functions, we obtain 
 $$
 (2\pi)^{-2} \|{\cal A}_L f\|_2^2 = \frac{1}{2}
 \sum_{\alpha \in \zd_+ }\left|\sum_{n=1}^d L_n \left(\sqrt{\alpha_n}c_{\alpha_1\dots\alpha_n-1\dots\alpha_d}   +\sqrt{\alpha_{n}+1}c_{\alpha_1\dots\alpha_n+1\dots\alpha_d}  \right)\right|^2.
 $$

Using the property of Hermite functions $\widehat{h_{\alpha}}(\xi)=(-\mathrm{i})^{|\alpha|} (2\pi)^{d/2} h_{\alpha}(2\pi \xi)$, we   
  get analogously 
\begin{eqnarray*}
  (2\pi)^2 \|{\cal B}_L f\|_2^2
 &=& (2\pi)^2
   \int_{\rd}\Bigl|\sum_{\alpha \in \zd_+ }c_{\alpha}\langle L, x \widehat{h_{\alpha}}\rangle\Bigr|^2 \, \mathrm{d}x = 
 \int_{\rd}\Bigl| \sum_{\alpha \in \zd_+ }(-\mathrm{i})^{|\alpha|} c_{\alpha}\langle L, x h_{\alpha}\rangle \Bigr|^2 \, \mathrm{d}x
\\
 & = &
 \frac{1}{2} 
 \int_{\rd}\Bigl|  \sum_{\alpha \in \zd_+ } (-\mathrm{i})^{|\alpha|-1} h_{\alpha}(x) 
  \sum_{n=1}^d L_n \left(\sqrt{\alpha_n}c_{\alpha_1\dots\alpha_n-1\dots\alpha_d}   -\sqrt{\alpha_{n}+1}c_{\alpha_1\dots\alpha_n+1\dots\alpha_d}
  \right)\Bigr|^2 \, \mathrm{d}x 
\\
& =& \frac 12 \sum_{\alpha \in \zd_+ }\left|\sum_{n=1}^d L_n \left(\sqrt{\alpha_n}c_{\alpha_1\dots\alpha_n-1\dots\alpha_d}    - \sqrt{\alpha_{n}+1}c_{\alpha_1\dots\alpha_n+1\dots\alpha_d}  \right)\right|^2.
\end{eqnarray*}
 Now, formula (\ref{hermite1}) immediately follows from above expressions for $(2\pi)^{-2} \|{\cal A}_L f\|_2^2$ and $(2\pi)^2 \|{\cal B}_L f\|_2^2$. 
To get (\ref{hermite}) we write 
 $$
(2\pi)^{-2} \|{\cal A}_L f\|_2^2 + (2\pi)^2 \|{\cal B}_L f\|_2^2 = 
 \frac{1}{2}
 \sum_{\alpha \in \zd_+ }\left(\left|\sum_{n=1}^d L_n \left(\sqrt{\alpha_n}c_{\alpha_1\dots\alpha_n-1\dots\alpha_d}   +\sqrt{\alpha_{n}+1}c_{\alpha_1\dots\alpha_n+1\dots\alpha_d}  \right)\right|^2\right.
 $$
 $$
 +\left.
 \left|\sum_{n=1}^d L_n \left(\sqrt{\alpha_n}c_{\alpha_1\dots\alpha_n-1\dots\alpha_d}   -\sqrt{\alpha_{n}+1}c_{\alpha_1\dots\alpha_n+1\dots\alpha_d}  \right)\right|^2
 \right)
 $$
 $$
 =
 \sum_{\alpha \in \zd_+ }
 \left(\left|\sum_{n=1}^d L_n \sqrt{\alpha_n}c_{\alpha_1\dots\alpha_n-1\dots\alpha_d} \right|^2+
 \left|\sum_{n=1}^d L_n \sqrt{\alpha_{n}+1}c_{\alpha_1\dots\alpha_n+1\dots\alpha_d}  \right|^2
 \right).
 $$
\hfill $\Box$

{\bf Proof of Lemma \ref{lem:symm}.}
 The symmetry relations (\ref{pm_symm}) mean that 
 $c_{\alpha_1\dots\alpha_n-1\dots\alpha_d} c_{\alpha_1\dots\alpha_k-1\dots\alpha_d}=0$ for 
 $k\neq n,$
 $k, n = 1,\dots,d$.
 
  So,
(\ref{hermite}) is rewritten as
$$
(2\pi)^{-2} \|{\cal A}_L f\|_2^2 + (2 \pi)^2 \|{\cal B}_L f\|_2^2   
 =
 \sum_{\alpha \in \zd_+ }
 \left(\sum_{n=1}^d L_n^2 \alpha_n \left|c_{\alpha_1\dots\alpha_n-1\dots\alpha_d} \right|^2+
 \sum_{n=1}^d L_n^2 (\alpha_{n}+1) \left|c_{\alpha_1\dots\alpha_n+1\dots\alpha_d}  \right|^2
 \right)
$$
$$
 =
 \sum_{\alpha \in \zd_+ }
 \sum_{n=1}^d L_n^2 (2 \alpha_n + 1 ) \left|c_{\alpha} \right|^2 \ge 
 3 \sum_{\alpha \in \zd_+ }
 \sum_{n=1}^d L_n^2  \left|c_{\alpha} \right|^2 = 3
 \|L\|^2 \|f\|_2^2.
$$
Since the function $c f(c\cdot)$ keeps the symmetry, it follows that 
$$
\|{\cal A}_L f\|_2 \|{\cal B}_L f\|_2 \geq \frac{3}{2}\|L\|^2 \|f\|_2^2.
$$
Finally, by (\ref{0center}), we obtain 
$
\displaystyle
\mathrm{UP}_L(f)\ge \frac{9}{4}.
$
\hfill $\Box$

 \section{Examples}

We give a couple of examples to illustrate the results of  Subsection 3.2, namely, the dependence of  localization on a direction $L$.

\textit{Example 1.} We illustrate Remark \ref{indep}. Let $d=2$, $L=(a,b)$, where $a,b\in\mathbb R$ and $a^2+b^2=1$. Consider the function 
$$
f(x,y):=\frac{3 xy}{2}\, \mathbbm{1}_{[-1,\,1]^2} (x,\,y).
$$
Note that $\|f\|_2=1$, 
$$
{\cal A}_L f = 3\pi (ax+by)xy,\ \ \ 
{\cal B}_L f = \frac{3\pi(ay+bx)}{4}.
$$
Since $\alpha_L(f)=\beta_L(f)=0$, it follows that 
\begin{equation*}
\begin{aligned}
(2\pi)^{-2}\Delta({\cal A}_L,f) + (2\pi)^2\Delta({\cal B}_L,f) &= (2\pi)^{-2} \|{\cal A}_L f\|_2^2 + (2\pi)^2 \|{\cal B}_L f\|_2^2 \\
&= 3(a^2+b^2)/5 + 3\cdot(a^2+b^2)= 18/5.
\end{aligned}
\end{equation*}
So $(2\pi)^{-2}\|\Delta({\cal A}_L,f)\|^2+ (2\pi)^2 \|\Delta({\cal B}_L,f)\|^2$ is constant and does not depend on $L$.

\textit{Example 2.} Let $d=2$,  $L=(a,b)$, where $a,b\in\mathbb R$ and $a^2+b^2=1$. Consider the function 
$$
f(x,y):=\frac{\sqrt{21} x^3y}{2} \,\mathbbm{1}_{[-1,\,1]^2} (x,\,y).  
$$
Note that $\|f\|_2=1$,
$$
{\cal A}_L f = \sqrt{21}\pi\, (ax+by)x^3y,\ \ \ 
{\cal B}_L f = \frac{\sqrt{21}}{4\pi}\,(3ax^2y+bx^3).
$$
Since $\alpha_L(f)=\beta_L(f)=0$,  it follows that 
\begin{equation*}
\begin{aligned}
(2\pi)^{-2} \Delta({\cal A}_L,f) + (2\pi)^2 \Delta({\cal B}_L,f) &= (2\pi)^{-2}\|{\cal A}_L f\|_2^2 + (2\pi)^2 \|{\cal B}_L f\|_2^2 \\
&= \left(\frac{7}{9}\, a^2 + \frac{3}{5}\, b^2\right) + \left(\frac{28}{5}\, a^2 + \frac{4}{3}\, b^2\right) = 
\frac{287}{45}\, a^2 + \frac{29}{15}\, b^2  \\
&=\frac{40}{9}\, a^2 + \frac{29}{15}.
\end{aligned}
\end{equation*}
So, $(2\pi)^{-2}\Delta({\cal A}_L,f) + (2\pi)^2 \Delta({\cal B}_L,f)$ is a quadratic function of $a$ and its maximum and minimum is attained on $a=1$ ($b=0$) and $a=0$ ($b=1$) respectively.

\textit{Example 3.} We illustrate Theorem \ref{depend_L}.  Consider a function
$$
f_0(x):=(2/\pi)^{d/4} (a_1\cdot...\cdot a_d)^{1/4} {\mathrm e}^{-(a_1 x_1^2+...+a_d x_d^2)}
$$
where $x \in\mathbb{R}^d$ and $a_1,...,a_d>0$. Let $L\in\rd$,  and $\|L\|=1$.
The Fourier transform of the function $f_0$ is
$$
\widehat f_0(x_1,...,x_d) = (2\pi)^{d/4}(a_1 \cdot...\cdot a_d)^{-1/4}{\mathrm e}^{-\pi^2(x_1^2/a_1 + ... + x_d^2/a_d)}.
$$
Since $f_0$ is even with respect to every variable $x_k,$ $k=1,\dots,d$, it follows that $M_{n,k}=0$ for $k\neq n$, $n, k=1,\dots,d$.

So, the matrix $M$ is diagonal and for the $k$-th diagonal element we obtain 
\begin{equation*}
\begin{aligned}
M_{k,k} &= \int_{\mathbb R^d}x_k^2|f(x)|^2\,\mathrm{d}x + 
(2\pi)^2\int_{\mathbb R^d}x_k^2|\widehat f(x)|^2\,\mathrm{d}x\\
&=
(2/\pi)^{d/2}(a_1 \dots a_d)^{1/2}\int_{\mathbb R^d}x_k^2 {\mathrm e}^{-2(a_1 x_1^2+...+a_d x_d^2)}\,\mathrm{d}x\\ 
&\ \ \ +
(2\pi)^{d/2+2}(a_1 \dots a_d)^{-1/2}\int_{\mathbb R^d}x_k^2 {\mathrm e}^{-2\pi^2( x_1^2/a_1+...+x_d^2/a_d)}\,\mathrm{d}x\\
&=\frac{1}{4a_k}+a_k.
\end{aligned}
\end{equation*}
Therefore, eigenvectors of $M$ coincide with the standard basis $\{e_k\}_{k=1}^d$ and corresponding eigenvalues are equal to $M_{k,k}$.
Thus, by Theorem \ref{depend_L} 
$$
\min_{\|L\|=1} \left((2\pi)^{-2}\Delta({\cal A}_L, f) + (2\pi)^2 \Delta({\cal B}_L, f)\right) = \min_{k=1,\dots,d} M_{k,k}, 
$$
$$
\max_{\|L\|=1} \left((2\pi)^{-2}\Delta({\cal A}_L, f) + (2\pi)^2\Delta({\cal B}_L, f)\right) = \max_{k=1,\dots,d} M_{k,k}.
$$

The  case of a function 
$$
f(x):=(2/\pi)^{d/4} (a_1 \dots  a_d)^{1/4}
{\mathrm e}^{-\sum_{i=1}^d\sum_{j=1}^d a_{i,j}x_i x_j},
$$
$a_{i,j}>0$, $i,j=1\dots,d$, reduces to the case of the function $f_0$. Indeed,   by suitable  shifts and   rotations the function $f$ transforms to $f_0$. Then by Lemma \ref{lem:UCRdSTD}, we conclude that these transformations do not change the maximum and the minimum of  $(2\pi)^{-2}\Delta(A_L,f) + (2\pi)^2 \Delta(B_L,f)$.

{\it Example 4.} We illustrate Theorem \ref{1depend_L}. Consider the same  function $f_0$ as in Example 3.  
The moments are
 $$
  M_k=(2\pi)^2 \int_{\rd}x_k^2 |f(x)|\,\mathrm{d}x = \frac{\pi^2}{a_k},
 \quad
  \widehat{M}_k=\int_{\rd}\xi_k^2 |\widehat{f}(\xi)|\,\mathrm{d}\xi = \frac{a_k}{4 \pi^2}.
  $$
  So, the elements of the matrix $A$ are equal to 
  $
  A_{n,k} = 1/8 \left(a_k/a_n+ a_n/a_k\right).$
  
  It turns out that determinants of the matrix $A$ and all the matrices $A_{j_1,\dots,j_q}$ are equal to zero for $q=1,\dots,d-2.$  
According to Theorem \ref{1depend_L}, it means that the set of all extremal vectors $v$ consists of the vectors with at least $d-2$ nonzero coordinates. So, these are vectors  of the type  
$v^{nk}:=(0,\dots,0,v_k,0\dots,0,v_n,0,\dots,0),$ $n,k=1,\dots,d$, where $v_k\neq 0$, $v_n\neq 0$,  and vectors $e_k$, $k=1,\dots,d$, of the standard basis in $\rd$. The vectors $v^{nk}$ satisfy the equation
$$
\frac18
\left(
\begin{array}{cc}
  2   & a_k/a_n + a_n/a_k\\
    a_k/a_n + a_n/a_k & 2
\end{array}
\right) 
\left(
\begin{array}{c}
  v_k\\
 v_n
\end{array}
\right) = \lambda 
\left(
\begin{array}{c}
  1\\
 1
\end{array}
\right).
$$
Since the parameter $\lambda$ is chosen to satisfy the condition  $v_1+\dots+v_d=1$, it  follows that $v_k=v_n=1/2.$ Therefore, the corresponding extremal directional vector $L^{nk}$ has two nonzero coordinates $L_n=L_k=1/\sqrt{2}.$ 
Calculating and comparing the values of $\mathrm{UP}_L(f)$ for the vectors $L^{nk}$ and $e_k$ 
we obtain 
$$
{\rm max}_{\|L\|=1} \mathrm{UP}_L(f) = {\rm max}_{n,k=1\dots,d} \mathrm{UP}_{L^{nk}}(f) = \frac{1}{16}
{\rm max}_{n,k=1\dots,d} \frac{(a_n+a_k)^2}{4 a_n a_k},
$$
$$
{\rm min}_{\|L\|=1} \mathrm{UP}_L(f) = \mathrm{UP}_{e_{k}}(f) = \frac{1}{4}.
$$



\bibliographystyle{elsarticle-num}


\end{document}